\documentclass[12pt]{article}
\usepackage{natbib}
\usepackage{hyperref}
\usepackage{amssymb}
\usepackage{dsfont}
\usepackage{hhline}
\usepackage{epsfig}
\begin{document}
\newtheorem{theorem}{Theorem}
\newtheorem{remarque}{Remark}
\newtheorem{lemma}{Lemma}
\newtheorem{corollaire}{Corollary}

%
%
\vskip 3mm

\vskip 3mm \baselineskip=15pt

\noindent {\Large Asymptotic Normality of the Additive Regression
Components for Continuous Time Processes}

\vskip 5mm \noindent Mohammed DEBBARH   and Bertrand MAILLOT\vskip
2mm \noindent

\noindent Universit\'e Paris 6
\\
\noindent 175,\ Rue du Chevaleret,  75013 Paris.\\
\noindent debbarh@ccr.jussieu.fr and maillot@ccr.jussieu.fr.
\\

\noindent ABSTRACT

In multivariate regression estimation, the rate of convergence
depends on the dimension of the regressor. This fact, known as the
{\it curse of the dimensionality}, motivated several works. The
additive model, introduced by Stone \cite{Stone1985}, offers an
efficient response to this problem. In the setting of continuous
time processes, using the marginal integration method, we obtain
the quadratic convergence rate and the asymptotic normality of the
components of the additive model.

\section{Introduction}
 Let ${\bf Z}_t=({\bf X}_t, Y_t)_{(t\in \mathbb{R})}$ be a
$\mathbb{R}^d\times \mathbb{R}$-valued measurable stochastic
process defined on a probability space $(\Omega, \mathcal{A}, P)$
with $d \geq 1$. Let $\mathcal{C}_1, ...,~ \mathcal{C}_d,$ be $d$
compact intervals of $\mathbb{R}$ and set
$\mathcal{C}=\mathcal{C}_1\times...\times \mathcal{C}_d$. Set now
$\delta >0$ and introduce the $\delta$-neighborhood
$\mathcal{C}^\delta$ of $\mathcal{C}$, namely $\mathcal{C}^\delta
= \{{\bf x} : \inf_{{\bf y}\in \mathcal{C}}\|{\bf x}-{\bf
y}\|_{\mathbb{R}^d}< \delta\},$ with $\|\cdot\|_{\mathbb{R}^d}$
standing for the euclidian norm on $\mathbb{R}^d$. Let $\psi$ be a
real valued measurable function. Consider the regression function
$m_\psi$ defined by, \vspace{-.2cm}
\begin{eqnarray}
m_{\psi}(\bf{x})& = & E\left(\psi (Y) \mid \bf{X}=\bf{x}\right),~
\forall ~{\bf x}=(x_1,...,x_d) \in \mathcal{C}^{\delta}.
\label{fdereg}
\end{eqnarray}
\vskip5pt\noindent Let $K$  be a kernel  defined on $\mathbb{R}^d$
and having a compact support. Let $\hat f_T$ be the estimate of
$f$, the density function of the covariable ${\bf X }$, (see Banon
\cite{Banon1978}), defined by,

\begin{eqnarray*}
\hat f_T({\bf x})= \frac{1}{Th_T^d} \int_0^T K\Big(\frac{{\bf
x}-{\bf X}_s}{h_T}\Big) ds,
\end{eqnarray*}
where $h_T$ is a given real positive function. In the sequel, to
estimate the regression function defined in (\ref{fdereg}), we use
the following  estimator (see, for example, Bosq \cite{Bosq1993}
and Jones et al.\cite{Jones1994})

\begin{equation}
\widetilde{m}_{\psi,T}({\bf x}) = \int_0^T W_{T,t}({\bf
x})\psi({\bf Y}_t)dt ~~ \mbox{with}~~W_{T,t}({\bf x})= \frac{
\prod_{l=1}^d\frac{1}{h_{l,T}}K_{l}\big({\frac{{\bf x}_{t} -{\bf
X}_{t}}{h_{l,T}}}\big)}{T\hat f_T ({\bf X}_t)},\label{estcarbo1}
\end{equation}
where $(h_{j,T})_{1\leq j \leq d}$ are positive real functions and
$(K_l)_{1\leq j \leq d}$ are $d$ kernels defined on $\mathbb{R}$
with compact supports. Consider now that the nonparametric
regression function (\ref{fdereg}) may be written as a sum of
univariate functions, i.e.

\begin{equation}\label{additive}
 m_\psi({\bf x}) \equiv  \mu + \sum_{l=1}^{d} m_{l}(x_{l})=: m_{\psi,add}({\bf x}),~\forall ~{\bf x}=(x_1,...,x_d) \in
 \mathcal{C}^{\delta},
\end{equation}
where, for $1\leq l \leq d$, $Em_l(X_l)=0$. For $1\leq l \leq d$
and any ${\bf x}=(x_1,..,x_{d})\in\mathcal{C}^{\delta}$ set ${\bf
x}_{-l}=(x_1,..,$ $x_{l-1},x_{l+1},$ $..,$ $x_d)$. To estimate the
additive components, we use the marginal integration method (see
Linton \& Nielsen \cite{Linton1995} and Newey \cite{newey1994}).
To this aim, we introduce $d$ densities $q_1,...,q_d$ defined on
$\mathbb{R}$ and set $q({\bf x}) = \prod_{l=1}^d q_l(x_l)$ and
$q_{-l}({\bf x}_{-l}) = \prod_{j \neq l} q_j(x_j)~\{l=1,...,d\}$.
We can then write
\begin{eqnarray}
m_{\psi}({\bf x})= \sum_{l=1}^d \eta_l(x_l) + \int_{\mathbb{R}^d}
m_{\psi}({\bf z})q({\bf z}) d{\bf z}
\label{additive_component_marginale}
\end{eqnarray}
with
\begin{eqnarray}
\mbox{with}~\eta_{l}(x_{l}) &:=& \int_{\mathbb{R}^{d-1}}
m_{\psi}({\bf x}) q_{-l}({\bf x}_{-l}) d{\bf x}_{-l} -
\int_{\mathbb{R}^d} m_{\psi}({\bf x}) q({\bf x}) d{\bf x}\nonumber \\
&=& m_l(x_l) - \int_{\mathbb{R}} m_l(z) q_l(z)dz, 1\leq l \leq d.
\label{relation_additive_component}
\end{eqnarray} Making use of the statements (\ref{estcarbo1}) and (\ref{relation_additive_component}), it follows that a natural
estimate of the $l$-th component is given by
\begin{eqnarray}
\widehat \eta_{l,T}(x_{l}) = \int_{\mathbb{R}^{d-1}}\!
\widetilde{m}_{\psi,T}({\bf x}) q_{-l}({\bf x}_{-l}) d{\bf x}_{-l}
- \int_{\mathbb{R}^d}\! \widetilde{m}_{\psi,T}({\bf x}) q({\bf x})
d{\bf x},\ l=1,...,d.
\end{eqnarray}
\vspace{-1.3cm}
\section{Hypotheses and Notations}
 \noindent In order to state our results, we introduce some
assumptions and additional notations.\\
\\
(C.1) $\mbox{There exists a positive constant }M \mbox{ such that,
for any
$y\in\mathbb{R},~~|\psi(y)| \leq M < \infty$},$\\
(C.2) $m_\psi \mbox{ is a $k$-times continuously differentiable
function, $ k \geq 1$, and }$\\  $\sup_{{\bf x}}
\Big|\frac{\partial^k m_\psi}{\partial x_l^k}({\bf x})
\Big|<\infty;~1\leq l \leq d.$\\
\\
For $1\leq l \leq d$, we denote by $f_l$,  the density function of
$X_l$ and we suppose that the functions $f$ and $f_l$ are
continuous and bounded. We need the additional conditions\\
\\
$(F.1)~ \forall{\bf x}\in\mathcal{C}^\delta, ~ f({\bf
x})>0~\mbox{and}~f_{l}(x_l)>0,\
l=1, ..., d,$  \\
$(F.2)~ f \mbox{ is $k'$-times continuously differentiable on
} \mathcal{C}^{\delta}, k'> k d,$  \\
$(F.3)~ \mbox{for some} ~~ 0 < \lambda \leq 1, \Big|\frac{\partial
f^{(k')}}{\partial x_1^{j_1}...\partial_d^{j_d}}({\bf
x}')-\frac{\partial f^{(k')}}{\partial
x_1^{j_1}...\partial_d^{j_d}}({\bf x})\Big| \leq L \|{\bf
x}'-{\bf x}\|^\lambda~\mbox{with}~j_1+...+j_d=k'.$\\
Here $\|.\|$ states as a norm on $\mathbb{R}^d$, $L$ is a positive
constant and we note $r:=k'+\lambda$.\\
\\
The kernels $K$ and $K_l, 1\leq l \leq d$  are assumed  to fulfill the following conditions\\
\\
(K.1)~  For $1\leq l \leq d$, $K$ and $K_l$  are continuous
 on compact supports $S$ and
$S_l\subset\mathcal{C}_l$, respectively, \\
(K.2)~  $\int K=1$ and $\int K_j=1, ~~1\leq l \leq d,$\\
(K.3)~  $\mbox{$\prod_{j=1}^d K_j$ is of order $k$,}$\\
(K.4)~  $\mbox{$K$ is of order $k'$.}$\\
\\
The known integration density functions  $q_l$, $1\leq l \leq d$,
satisfy
the following assumption\\
\\
$(Q.1)$~ $q_l$ has  k  continuous and bounded derivatives, with
compact support included in
$\mathcal{C}_l,~~  1\leq l \leq d.$\\
\\
There exists $\Gamma \in \mathcal{B}_{\mathbb{R}^2}$ containing
$D=\{(s,t)\in \mathbb{R}^2: s=t\}$ such that\\
\\
$(D.1)~  f_{({\bf X}_s,Y_s),({\bf X}_t,Y_t)} -
f_{({\bf X}_s,Y_s)}\bigotimes f_{({\bf X}_t,Y_t)}~~\mbox{exists everywhere for} ~~ (s,t) \in \Gamma^C,$ \\
$(D.2)~  A_\Gamma:= \sup_{(s,t) \in \Gamma^C}\sup_{{\bf
x,y}\in\mathcal{C^{\delta}}\times\mathcal{C^{\delta}}}\int_{u,v\in\mathbb{R}^2}
| f_{({\bf X}_s,Y_s),({\bf X}_t,Y_t)}({\bf x},u,{\bf y},v) -$\\
~~~~~~~~~~~~~~~~~~~~~~~~~~~~~~~~~~~~~~~~$f_{({\bf X}_s,Y_s)}({\bf
x},u)
f_{({\bf X}_t,Y_t)}({\bf y},v)|dudv <\infty,$\\
$(D.3)~\mbox{there exists} ~\ell_\Gamma<\infty~\mbox{and} ~T_0~
\mbox{such that},~ \forall T>T_0,~ \frac{1}{T}
\int_{[0,T]^2\cap \Gamma} ds dt\leq \ell_\Gamma.$ \\
\\
We will work under the following conditions on the smoothing
parameters $h_T$ and $h_{j,T},~j=1,...,d$.\\
\\
$(H.1)~ h_T =c' \Big(\frac{\log T}{T}\Big)^{1/(2k'+d)}, \mbox{ for a fixed } 0<c'<\infty$,\\
$(H.2)~ h_{j,T}=c_1T^{-1/(2k+1)}, \mbox{ for fixed } 0<c_1<\infty.$\\
\\
Let $\mathcal{A}$  and $\mathcal{B}$ be two $\sigma$-fields. We
will use the $\alpha$-mixing coefficient defined by
$$\alpha\big(\mathcal{A},\mathcal{B}\big)=\sup_{(A,B)\in(\mathcal{A},\mathcal{B})}|P(A\cap
B)-P(A)P(B)|.$$ For all Borel set $I\subset\mathbb{R}^+$ the
$\sigma$-algebra defined by $\big(Z_t,t\in I\big)$ will be denoted
by $\sigma\big(Z_t,t\in I\big)$. Writing
$\alpha(u)=\sup_{t\in\mathbb{R}_+}\alpha\big(\sigma\big(Z_v,v\leq
t\big),\sigma\big(Z_v,v\geq t+u\big)\big)$, we will use the
condition \\
\\
$(A.1)~\alpha(t)=\mathcal{O}\big(t^{-b}\big) ~~ \mbox{with}~~b >
\frac{7r+5d}{2r}$.\\
\\
We denote by $\widehat{\widehat{\eta}}_{l,T}$ and
$\widetilde{\widetilde{m}}_{\psi,T}({\bf x})$ the versions of
$\widehat{\eta}_{l,T}$ and $\widetilde{m}_{\psi,T}({\bf x})$
corresponding to a known density $f$. Introduce now the following
quantities (see, for the discrete case, Camlong {\it et al.}
\cite{Vieu2000}),

\begin{eqnarray*}
&& \hspace{-1cm} \tilde Y_{\psi,T,t,l}= \psi (Y_t)
\int_{\mathbb{R}^{d-1}} \prod_{j\neq l}^d \frac{1}{h_{j,T}}
K_{j}\Big({\frac{{\bf x}_{j} - {\bf X}_{t,j}}{h_{j,T}}}\Big)
\frac{q_{-l}({\bf x}_{-l})}{f(X_{t,-l}|X_{t,l})} d{\bf x}_{-l};
\widetilde m_{\psi,l}^T(x_l)= E\big(  \tilde
Y_{\psi,T,t,l}\Big| X_{t,l} = x_{l}\big);\\
&&\hspace{-1cm} \widehat \alpha_{l}(x_{l}) = \frac{1}{Th_{l,T}}
\int_0^T \frac{\tilde Y_{\psi,T,t}}{f_1(X_{t,l})} K_{l}
\Big(\frac{x_{l}- X_{t,l}}{h_{l,T}}\Big)dt; \mathcal{G}_{l}({\bf
u}_{-l})=\int_{\mathbb{R}^{d-1}} \prod_{j\neq l}^d
\frac{1}{h_{j,T}} K_{j}\Big({\frac{{\bf x}_{j} - {\bf
u}_{j}}{h_{j,T}}}\Big)
q_{-l}({\bf x}_{-l}) d{\bf x}_{-l};\\
&&\hspace{-1cm} C_{T,l} =  \mu + \int_{\mathbb{R}^{d-1}}
\sum_{j\neq l} m_j(u_j) \mathcal{G}_{l}({\bf u}_{-l}) d{\bf
u}_{-l}; \widehat C_T=\int_{\mathbb{R}^d} \widetilde{\widetilde
m}_{\psi,T}({\bf x})
q({\bf x})d{\bf x}; C_{l}=\int_{\mathbb{R}}m_{l}(x_{l})q_{l}(x_{l})dx_{l};\\
&& b_l(x_l) = \frac{1}{k!}\int_{\mathbb{R}}u^kK_l(u)du \Big((-1)^k
m_l^{(k)}(x_l)+\int_{\mathbb{R}}m_l(z)q_l^{(k)}(z)dz \Big).
\label{def_b_l}
\end{eqnarray*}

\section{Results}
The proofs of our Theorems are split into two steps. We first
consider the density as known, and then treat the general case
where $f$ is unknown by using the  decomposition $1/f= 1/\hat f_T
- (f-\hat f_T)/f\hat f_T$ and the following lemma.
\begin{lemma}\label{lem} Under the assumptions $(F.1)-(F.3)$,
$(K.1),(K.2),(K.4)$, $(D.1)-(D.3)$, $(H.1)$ and $(A.1)$we have
\begin{eqnarray}\label{bosq_result}
\sup_{{\bf x} \in \mathcal{C}} |\hat f_T ({\bf x}) - f({\bf x})|=
\mathcal{O}\Big(\Big(\frac{\log
T}{T}\Big)^{k'/(2k'+d)}\Big)~~\mbox{a.s.}.
\end{eqnarray}
\end{lemma}
{\bf Proof:} It is easily seen that under our assumptions, the
result follows by using the arguments used in the  demonstration
of Theorem 4.9. in \cite{Bosq1996} p.112 and by replacing $\log_m$
by $1$.
\begin{theorem}\label{lem1}
Under assumptions $(C.1)-(C.2)$, $(F.1)-(F.3)$, $(K.1)-(K.4)$,
$(Q.1)$,  $(D.1)-(D.3)$, $(H.1)-(H.2)$ and $(A.1)$ we have
\begin{eqnarray*}
 & &E\big( {\widehat \eta}_{l,T}(x_l) -
\eta_{l}(x_l)\big)^2 = \mathcal{O}\Big(T^{-2k/(2k+1)}\Big).
\label{equation1}
\end{eqnarray*}
\end{theorem}
{\bf Sketch of the proof:} Observe that \vspace{-.2cm}
\begin{eqnarray}
{\widehat \eta}_{l,T}(x_l) - \eta_l(x_l) & = & \{{\widehat
\eta}_{l,T}(x_l)-\widehat{\widehat \eta}_{l,T}(x_l)\}+ \{\hat
\alpha_l(x_l) - E\hat \alpha_l(x_l)\} +\{E\hat \alpha_l(x_l)-
\tilde m_{\psi,l}^T(x_l)\}\\&& +E\{\hat C_T - C_{T,l} -C_l\}.
\nonumber \label{decomp_biais_eta}
\end{eqnarray}
It follows that
\begin{eqnarray*}
\hspace{-.8cm}E\{{\widehat \eta}_{l,T}(x_l) - \eta_l(x_l)\}^2 &
\leq & 4E\{{\widehat \eta}_{l,T}(x_l)-\widehat{\widehat
\eta}_{l,T}(x_l)\}^2+ 4 E\{\hat \alpha_l(x_l) - E\hat
\alpha_l(x_l)\}^2+4\{E\hat \alpha_l(x_l)- \tilde
m_{\psi,l}^T(x_l)\}^2\\&&+4E^2\{\hat C_T - C_{T,l} -C_l\}.
\end{eqnarray*}
To prove  the Theorem \ref{lem1}, it suffices to establish the
following statements \vspace{-.2cm}
\begin{eqnarray}
&&E({\widehat \eta}_{l,T}(x_l)-\widehat{\widehat
\eta}_{l,T}(x_l))^2
=\mathcal{O}\Big(T^{-2k/(2k+1)}\Big), \label{decom_1}\\
&&{\rm Var}(\hat
\alpha_l(x_l))=\mathcal{O}\Big(T^{-2k/(2k+1)}\Big),
\label{equation3}~\\
&&E\hat \alpha_l(x_l)- \tilde
m_{\psi,l}^T(x_l)=\mathcal{O}\Big(T^{-k/(2k+1)}\Big)\label{decomp_lem1_3},\\
&&E(\hat C_T -C_{T,l}+C_l)
=\mathcal{O}\Big(T^{-k/(2k+1)}\Big).\label{equation2}
\end{eqnarray}
{\it Proof of \ref{decom_1}:} By combining the definitions of
${\widehat \eta}_{1,T}$ and $\widehat {\widehat \eta}_{1,T}$ and
the result of the lemma \ref{lem}, we easily obtain, under the
conditions on the kernel, the statement(\ref{decom_1}).\\
{\it Proof of \ref{equation3}:} Set $\phi(t,s)= {\rm
Cov}\Big(\frac{\widetilde
Y_{\psi,T,t}}{f_{1}(X_{t,1})h_{1,T}}K_1\Big(\frac{x_1-X_{t,1}}{h_{1,T}}\Big),
 \frac{\widetilde
Y_{\psi,T,s}}{f_{1}(X_{s,1})h_{1,T}}K_1\Big(\frac{x_1-X_{s,1}}{h_{1,T}}\Big)\Big)
$ and $S_{a(T)}=\{(s,t)\in\mathbb{R}^2;|t-s|\leq a(T)\}$, where
$a(T) = h_T^{-1}$. We use the following decomposition
\vspace{-.2cm}
\begin{eqnarray*}
\hspace{-.8cm}{\rm Var} (\hat \alpha_1(x_1))&=& \int_{[0,T]^2\cap
\Gamma}\phi(t,s)dtds +\int_{[0,T]^2\cap \Gamma^c \cap
S_{a(T)}}\phi(t,s)dtds+\int_{[0,T]^2\cap \Gamma^c \cap
S_{a(T)}^c}\phi(t,s)dtds:=A+E+F.
\end{eqnarray*}
Under $(C.1)$, $(F.1)$, $(K.1)-(K.2)$ and $(Q.1)$, we have, for
$T$ large enough, \vspace{-.2cm}
\begin{eqnarray}\label{AE}
A= \mathcal{O}\Big(1/Th_{1,T}\Big) ~~\mbox{and}~~E=
\mathcal{O}\Big(a(T)\|K_1\|_{\mathbb{L}_1}^2A_f(\Gamma)/T\Big).
\end{eqnarray}
Using the Billingsley's inequality, it follows that \vspace{-.2cm}
\begin{eqnarray}\label{fff}
F=\mathcal{O}\Big(1/Th_{1,T}^{2}a(T)\Big).
\end{eqnarray}
Combining (\ref{AE}) and (\ref{fff}), we obtain (\ref{equation3}).
To prove the statements  (\ref{decomp_lem1_3}) and
(\ref{equation2}), we use  similar arguments  as in  the discrete
case (see Camlong {\it
et al.} \cite{Vieu2000}).\\
\\
The next Theorem needs the following additional hypothesis.
\\ $(\mathcal{V})$ ~~~~~~~$\liminf_{T \rightarrow \infty } Th_{l,T}{\rm Var}( \hat
\eta_{l,T}(x_l))>0$ where $(\log (T)
/T)^{k'/(2k'+d)}=o(h_{l,T}^{k})$.
\begin{theorem}\label{th}
Under the hypotheses of Theorem \ref{lem1} and $(\mathcal{V})$ we
have, for every $\forall l\in [1, d]$  and $\forall
x_l\in\mathcal{C}_l$,
\begin{eqnarray*}
\frac{\widehat \eta_{l,T} (x_l) - \eta_l(x_l)-
h_{l,T}^{k}b_l(x_l)}{\sqrt{{\rm Var} (\hat \eta_{l,T} (x_l))}}
\stackrel{\mathcal{L}}{ \longrightarrow} \mathcal{N}(0,1).
\end{eqnarray*}
\end{theorem}
{\bf Sketch of the proof:} To obtain our theorem it suffices to
show that \vspace{-.2cm}
\begin{eqnarray}
&& \sup_{x_l \in \mathcal{C}_l}|{\widehat
\eta}_{l,T}(x_l)-\widehat{\widehat
\eta}_{l,T}(x_l)|= \mathcal{O}\Big(\sup_{\bf x\in \mathcal{C}}|\hat f_T({\bf x}) - f({\bf x})|\Big) ~~\mbox{a.s.},\label{equart_estimateur}\\
&&\frac{\{\hat \alpha_l (x_l) - E(\hat \alpha_l
(x_l))\}}{\sqrt{{\rm Var}(\hat \alpha_l (x_l))}} \longrightarrow
\mathcal{N}\big(0, 1\big), \label{norm1}\\
&& E\hat \alpha_l(x_l)- \tilde m_{\psi,l}^T(x_l)=
\frac{(-h_{l,T})^k}{k!} m_l^{(k)}(x_l)\int_{\mathbb{R}} v_l^k
K_l(v_l)dv_l +
o(h_{l,T}^k),\label{equ5}\\
&&\mbox{and~~} E\{\hat C_T -C_{T,l}+C_l\}=\frac{h_{l,T}^k}{k!}
\int_{\mathbb{R}}q_l^{(k)}(x_l) m_l(x_l)dx_l
\int_{\mathbb{R}}v_l^kK_l(v_l)dv_l+ o(h_{l,T}^k). \label{equ6}
\end{eqnarray}
{\it Proof of \ref{equart_estimateur}:} The result arises directly from the definitions of estimates of $\eta_l$ and the conditions on the kernels $K_l, 1\leq l \leq d$.  \\
{\it Proof of \ref{norm1}:} Set $\frac{\{\hat \alpha_l (x_l) -
E(\hat \alpha_l (x_l))\}}{\sqrt{{\rm Var}(\hat \alpha_l (x_l))}} =
\int_0^T Z_t dt =: S_T$. We employ then the big block–-small block
procedure. Indeed setting, $ S_T = \sum_{j=1}^{k-1} (\nu_j+ \xi_j)
=: S_T' + S_T''$ where $\nu_j = \int_{j(p+q)}^{j(p+q)+p} Z_t dt$
and $ \xi_j=\int_{j(p+q)+p}^{(j+1)(p+q)} Z_t dt.$ Now, it suffices
to prove the following statements,
\begin{eqnarray}
&&ES_T''^2 \rightarrow 0~~\mbox{as}~~T \rightarrow +\infty, \label{equ1}\\
&&\Big|E(e^{itS_T'})-\prod_{j=0}^{k-1}E(e^{it\nu_j})\Big|\rightarrow
0~~\mbox{as}~~T \rightarrow +\infty,
\label{equ2}\\
&&\sum_{j=0}^{k-1}E[\nu_j^2]\rightarrow 1~~\mbox{as}~~T \rightarrow +\infty, \label{equ3}\\
&&\mbox{and~~}\sum_{j=0}^{k-1}E[\nu_j^2 \mathbb{I}_{\{\nu_j^2
>\epsilon \}}]\rightarrow 0~~\mbox{as}~~T \rightarrow +\infty.
\label{equ4}
\end{eqnarray}
To show  (\ref{equ3}) et (\ref{equ4}),  we use the same arguments
as those deployed in the discrete case.
\begin{lemma}\label{prop}
Under the conditions $(C.1)-(C.4)$, $(F.1)-(F.2)$, $(K.1)$,
$(Q.1)-(Q.2)$ and $(H.1)-(H.2)$, we have, for every $1\leq l \leq
d$ and for any $x_l\in\mathcal{C}_l$ and every
$(\alpha,\beta)\in]0;0,5[\times]0,5;1[$,
\begin{eqnarray}
\liminf_{T\rightarrow\infty} P\Big(T^{\frac{k}{2k+1}}\{\widehat
\eta_{l,T} (x_l) - \eta_l(x_l)- h_{l,T}^{k}b_l(x_l)\}
\in[Aq_{\alpha};Aq_{\beta}]\Big)\geq \beta-\alpha,
\end{eqnarray}
where $A:=\big(\limsup_{T\rightarrow
+\infty}T^{\frac{2k}{2k+1}}{\rm Var} (\hat \eta_{l,T}
(x_l))\big)^{1/2}$ and $q_u$ is such that
$P(\mathcal{N}(0,1)<q_{u})=u$.
\end{lemma}

\end{document}